\documentclass[twoside,10pt,leqno]{amsart}
\everymath={\displaystyle}
\usepackage{amsfonts}
\usepackage{amsmath}
\usepackage{amscd}
\usepackage{amssymb}
\usepackage{amsthm}
\usepackage{amsrefs}
\begin{document}
\newtheorem{theorem}{Theorem}
\newtheorem{lemma}{Lemma}
\newtheorem{corollary}{Corollary}
\newtheorem{conjecture}{Conjecture}
\newtheorem{prop}{Proposition}
\numberwithin{equation}{section}
\newcommand{\dif}{\mathrm{d}}
\newcommand{\intz}{\mathbb{Z}}
\newcommand{\ratq}{\mathbb{Q}}
\newcommand{\natn}{\mathbb{N}}
\newcommand{\comc}{\mathbb{C}}
\newcommand{\rear}{\mathbb{R}}
\newcommand{\prip}{\mathbb{P}}
\newcommand{\uph}{\mathbb{H}}

\title{\bf Large Sieve Inequality with Characters to Square Moduli}
\author{Liangyi Zhao}
\maketitle

\begin{abstract}
In this paper, we develop a large sieve type inequality with characters to square moduli.  One expects that the result should be weaker than the classical inequality, but, conjecturally at least, not by much.  The method is generalizable to higher power moduli.
\end{abstract}

\section{Introduction and historical background}

It was in 1941 that Ju. V. Linnik \cite{JVL1} originated the idea of large sieve, and he also made application to the distribution of quadratic non-residues.  A. R\'enyi studied the large sieve extensively and made important applications to the Goldbach problem.  Refinements in that direction have later been made by many. \newline

A set of real numbers $\{ x_k \}$ is said to be $\delta$-spaced modulo 1 if $x_j-x_k$ is at least $\delta$ away from any integer, for all $j \neq k$.  Hence the set $\{ x_k \}$ must be finite with cardinality not exceeding $\delta^{-1}$.  Throughout we assume that $0 < \delta \leq \frac{1}{2}$.  \newline

The large sieve inequality, which we henceforth refer to as the classical large sieve inequality, is stated as follows.  Different elegant proofs of the theorem can be found in \cite{HD}, \cite{PXG}, \cite{HM2}, \cite{HM}.  The theorem, in the following form, was first introduced by Davenport and Halberstam, \cite{DH1} and \cite{DH2}.

\begin{theorem}
Let $\{ a_n \}$ be an arbitrary sequence of complex numbers, $\{ x_k \}$ be a set of real numbers which is $\delta$-spaced modulo 1, and $M$, $N$ be integers with $N>0$.  Then we have
\begin{equation}
\sum_k \left| \sum_{n=M+1}^{M+N} a_n e \left( x_k n \right) \right|^2 \ll \left( \delta^{-1} +N \right) \sum_{n=M+1}^{M+N} |a_n|^2,
\end{equation}
where the implied constant is absolute.
\end{theorem}
Save for the more precise implied constant, the above inequality is the best possible.  Moreover, Cohen and Selberg have shown independently that 
\begin{equation} \label{cohenselb}
\sum_k \left| \sum_{n=M+1}^{M+N} a_n e \left( x_k n \right) \right|^2 \leq \left( \delta^{-1} -1+N \right) \sum_{n=M+1}^{M+N} |a_n|^2,
\end{equation}
which is absolutely the best possible, since Bombieri and Davenport \cite{BD1} gave examples of $\{x_k\}$ and $a_n$, with $\delta \to 0$, $N \to \infty$ and $N \delta \to \infty$ such that asymptotic equality holds in \eqref{cohenselb}. \newline

 This theorem admits corollaries for additive and multiplicative characters.  We derive
\begin{equation*}
\sum_{q=1}^Q \sum_{\substack{a \; \bmod{ \; q} \\ \gcd(a,q)=1}} \left| \sum_{n=M+1}^{M+N} a_n e \left( \frac{a}{q} n \right) \right|^2 \ll (Q^2+N) \sum_{n=M+1}^{M+N} |a_n|^2, 
\end{equation*}
and
\begin{equation*}
\sum_{q=1}^Q \frac{q}{\varphi(q)} \sideset{}{^{\star}}\sum_{\chi \; \bmod{ \; q}} \left| \sum_{n=M+1}^{M+N} a_n \chi (n) \right|^2 \ll (Q^2+N) \sum_{n=M+1}^{M+N} |a_n|^2,
\end{equation*}
where here and after, $\sideset{}{^{\star}}\sum$ means that the sum runs over primitive characters modulo the specified modulus only.  As usual, $\varphi(q)$ is the Euler $\varphi$ function. \newline

In this paper, we shall establish large sieve inequality for additive characters in which the moduli are squares, {\it id est}, $q^2$ rather than $q$, see Theorem~\ref{larsiev} in Section 4.  The problem reduces down to the spacing properties of rational numbers with square denominators.  The key idea that we employ in resolving such a problem is the Weyl-Hardy-Littlewood method for exponential sums.  It is also worthwhile to note that the method that we use in this paper may be generalized to higher power moduli.  However, as Weyl's estimates for exponential sums weaken when the polynomial in the amplitude is of high degree, our corresponding results are also weakened.  From Theorem~\ref{larsiev} we derive a corresponding large sieve type inequality for multiplicative primitive characters, see Corollary~\ref{larsievcor} in Section 4. \newline

The author wishes to thank Mr. Waldeck Sch\"{u}zter, a fellow graduate student at Rutgers University, for his help in writing the C++ program that generated the data of Table 1 in Section 6.  The author also thanks his thesis adviser, Henryk Iwaniec, who first suggested this problem to the author and who, of his advise and support, has been most generous.  The author thanks the referee for pointing out a mistake in an earlier version of the paper. \newline

The following notations and conventions are used throughout paper. \newline

\noindent $e(z) = \exp (2 \pi i z) = e^{2 \pi i z}$. \newline
$f = O(g)$ means $|f| \leq cg$ for some unspecified positive constant $c$. \newline
$f \ll g$ means $f=O(g)$. \newline
$f \asymp g$ means $c_1g \leq f \leq c_2 g$ for some unspecified positive
constants $c_1$ and $c_2$.  Unless otherwise stated, all implied constants in $\ll$, $O$ and $\asymp$ are absolute. \newline
$\| x \| = \inf \{ |x-k| \; : \; k \in \mathbb{Z} \}$ denotes the
distance of a real number $x$ to its closest integer. \newline
$\qedsymbol$ denotes the end of a proof or the proof is easy and standard.

\section{Heuristics and ``trivial'' bounds}

As stated earlier, we are interested in having an estimate of the following kind:
\begin{equation} \label{int}
\sum_{q=1}^Q \sum_{\substack{a  \; \bmod{ \; q^2} \\ \gcd(a,q)=1}} \left| \sum_{n=M+1}^{M+N} a_n e \left( \frac{a}{q^2} n \right) \right|^2 \ll \Delta \sum_{n=M+1}^{M+N} |a_n|^2.
\end{equation}
It is a remark attributed to Borel that the rational numbers in the real line are like stars in the heavens ``to illuminate the mystery of the continuum.''  Indeed, we must investigate the well-spacedness of some of these ``stars.''  As it certainly suffices to consider only $q$'s in dyadic intervals, we set
\begin{equation*}
S_Q= \left\{ \frac{a}{q^2} \in \ratq : \gcd (a, q)=1, \; 1 \leq a < q^2, \; Q < q \leq 2Q \right\}.
\end{equation*}
We easily see that if $x$ and $x'$ are two distinct elements of $S_Q$, then $\| x-x' \| \geq Q^{-4}$.  Therefore, just from the classical large sieve inequality, we may take 
\begin{equation} \label{Q4}
\Delta = Q^4+N
\end{equation}
in \eqref{int}.  On the other hand, all rational numbers $\frac{a}{q^2}$ in $S_Q$ with a fixed denominator $q^2$ are clearly $q^{-2}$-spaced.  Therefore, again by the virtue of the classical large sieve inequality, we may also take in \eqref{int} after summing over $q$,
\begin{equation}
\Delta=Q(Q^2+N). \label{Q3}
\end{equation}
It is worthwhile to note that when $N \asymp Q^3$, both \eqref{Q4} and \eqref{Q3} can be interpreted as $Q^4$.  However, neither \eqref{Q4} nor \eqref{Q3} exploit the fact that squares are so sparsely distributed among the integers.  One easily deduces that there are $\asymp Q^3$ rational numbers between 0 and 1 with square denominators and height at most $Q^2$.  Hence, these rational numbers are ``on average'' $Q^{-3}$-spaced.  Therefore, we aim to exploit these ``facts'' and improve the estimates in \eqref{Q4} and \eqref{Q3}.  Toward that end, we ``divide and conquer.''

\section{Preliminary Lemmas}

We begin by quoting the duality principle, which says that the norm of a bounded linear operator in a Banach space is the same as that of its adjoint operator.  To us, it amounts to the swapping of order of summations.  More precisely, we have

\begin{lemma}[Duality Principle]\label{dual}
Let $T=[t_{mn}]$ be a finite square matrix with entries from the complex numbers.  The following two statements are equivalent:
\begin{enumerate}
\item For any complex numbers $\{ a_n \}$, we have
\begin{equation*}
\sum_m \left| \sum_n a_n t_{mn} \right|^2 \leq D \sum_n |a_n|^2.
\end{equation*}
\item For any complex numbers $\{ b_n \}$, we have
\begin{equation*}
\sum_n \left| \sum_m b_m t_{mn} \right|^2 \leq D \sum_m |b_m|^2.
\end{equation*}
\end{enumerate}
\end{lemma}

\begin{proof}  This is quoted from \cite{HI4}. \end{proof}

We shall need the Poisson summation formula.  This asserts that if $f(x)$ is a reasonably well-behaved function, then summing $f(n)$ over all integers $n$ is the same as summing the Fourier transform of $f(x)$ over all integers $n$.  More precisely, we have

\begin{lemma} [Poisson Summation Formula] \label{poisum}
Let $f(x)$ be a function on the real numbers that is piece-wise continuous with only finitely many discontinuities and for all real number, $a$, satisfies
\begin{equation*}
f(a) = \frac{1}{2} \left[ \lim_{x \to a-} f(x) + \lim_{x \to a+} f(x) \right].
\end{equation*}
Moreover, $f(x) \ll (1+|x|)^{-c}$ for some $c>1$ with an absolute implied constant.  Then we have
\begin{equation*}
\sum_{n=-\infty}^{\infty} f(n) = \sum_{n=-\infty}^{\infty} \hat{f}(n),
\end{equation*}
where
\begin{equation*}
\hat{f}(x) = \int_{-\infty}^{\infty} f(y) e(xy) \dif y,
\end{equation*}
the Fourier transform of $f(x)$.
\end{lemma}
\begin{proof} This is quoted in as \cite{DB}. \end{proof}

Both of the above lemmas are proved using standard means.  We shall not succeed in proving our contention without the following lemma.
\begin{lemma}[Weyl Shift] \label{weylk}
Let $I$ be an interval of length $N$ and $f(x)$ be a polynomial of degree $k \geq 2$ with real coefficients.  Set $\kappa=2^{k-1}$ and let the leading coefficient, the coefficient of $x^k$, of $f(x)$ be $\alpha$.  Also set
\begin{equation*}
S= \sum_{n \in I} e(f(n)).
\end{equation*}
Then we have
\begin{equation*}
|S|^{\kappa} \leq 2^{2\kappa}
N^{\kappa-1} + 2^{\kappa} N^{\kappa-k} \sum_{r_1,\cdots,r_{k-1}} \min \left( N, \frac{1}{\| \alpha k! r_1 \cdots r_{k-1} \|} \right),
\end{equation*}
where each $r$ runs from 1 to $N-1$.
\end{lemma}
\begin{proof} This is Lemma 5.6 in \cite{ET}. \end{proof}

In short, if $f(x)$ is a linear polynomial, then $S$ is none other than a geometric series.  But if $f(x)$ is a polynomial of degree $k>1$, then $f(x) - f(x+h)$ will be a polynomial of degree $k-1$ in $x$.  If we iterate this process $k-1$ times, we get a geometric series and obtain some saving in the estimate of the modulus of $S$. \newline

\section{The main contention}

Again, we are interested in having an estimate of the following kind:
\begin{equation*}
\sum_{q=1}^Q \sum_{\substack{a  \; \bmod{ \; q^2} \\ \gcd(a,q)=1}} \left| \sum_{n=M+1}^{M+N} a_n e \left( \frac{a}{q^2} n \right) \right|^2 \ll \Delta \sum_{n=M+1}^{M+N} |a_n|^2.
\end{equation*}

Before we state and prove our main contention of this section, we first estimate 
\begin{equation} \label{del1}
M(Q,N) = \max_{x \in S_Q} \# \left\{ x' \in S_Q : \| x-x' \| < \frac{1}{2N} \right\},
\end{equation}
which is the central issue of our theorem.
\begin{lemma} \label{lsscent}
Given $\epsilon >0$ and $N \in \natn$, then we have
\begin{equation} \label{lsscenteq}
M(Q,N) \ll \frac{Q^3}{N}+ \left( \sqrt{Q} + \frac{Q^2}{\sqrt{N}} \right) N^{\epsilon},
\end{equation}
where the implied constant in \eqref{lsscenteq} depends on $\epsilon$ alone.
\end{lemma}
\begin{proof} The task before us is as follows.  Let $x=\frac{a}{q^2}$ and $x'=\frac{a_1}{q_1^2}$, with $\gcd(a,q)=\gcd(a_1,q_1)=1$.  Let $aq_1^2-a_1q^2 \equiv b \pmod{q^2q_1^2}$ with $|b| \leq \frac{1}{2} q^2q_1^2$.  \newline 

We have 
\begin{equation} \label{inq1}
0 \leq \left\| \frac{a}{q^2} - \frac{a_1}{q_1^2} \right\| = \frac{|b|}{q^2q_1^2}  < \frac{1}{2N}.
\end{equation}
This yields that $|b| \ll q^2 Q^2 N^{-1} =B$, say.  We want to estimate, for each $\frac{a}{q^2}$, the number of fractions $\frac{a_1}{q_1^2} \in S_Q$ satisfying \eqref{inq1}.  Let $b=aq_1^2-a_1q^2$.  We have $|b| <B $ and
\begin{equation} \label{ineq2}
\left\{ \begin{array}{cccc} q^2 & \equiv & -b \overline{a}_1, & \pmod{q_1^2} \\
b & \equiv & aq_1^2, & \pmod{q^2} \end{array}, \right.
\end{equation}
where $\overline{a}_1$ is the multiplicative inverse of $a_1$ modulo $q_1^2$. \newline

We shall estimate the number of $b$'s and $q_1$'s that satisfy the second congruence relation in \eqref{ineq2} with $q_1<Q$ and $|b| < B$, which clearly majorizes the maximum that we need to estimate in \eqref{del1}. \newline

First, we set $\phi (x) = \left( \frac{\sin \pi x}{2x} \right)^2$, a constant multiple of F\'ejer kernel.  We note that $\phi (x)$ is non-negative, $\phi(x) \geq 1$ for $|x| \leq \frac{1}{2}$ and $\phi (0) = \pi^2/4$.  Therefore
\begin{equation} \label{estm}
\sum_{n \equiv aq_1^2  \; \bmod{\; q^2}} \phi \left( \frac{n}{2B} \right)
\end{equation}
majorizes $M(Q)$.  There is certainly no unique choice for this test function $\phi(x)$.  However, as we shall presently apply Poisson summation formula, we find it most convenient to choose $\phi(x)$ this way, since its Fourier transform is a function of compact support, specifically $\hat{\phi}(s)=\frac{\pi^2}{4} \max(1-|s|,0)$. \newline

Now we apply Poisson summation, Lemma~\ref{poisum} with a linear change of variable, to \eqref{estm} and sum $q_1$ over dyadic intervals, we get
\begin{equation}
\frac{2B}{q^2} \sum_{Q < q_1 \leq 2Q} \sum_{j} e \left( \frac{ajq_1^2}{q^2} \right) \hat{\phi} \left( \frac{2jB}{q^2} \right).
\end{equation}

More precisely, the above is
\begin{eqnarray*}
 & \frac{\pi^2B}{2q^2} \sum_{|j| < \frac{q^2}{2B}} \sum_{Q < q_1 \leq 2Q} \left( 1 - \frac{2|j|B}{q^2} \right) e \left( \frac{ajq_1^2}{q^2} \right) \\
= & \frac{\pi^2 Q^2}{2N} \sum_{|j| < \frac{N}{4Q^2}} \sum_{Q < q_1 \leq 2Q} \left( 1 - \frac{4|j|Q^2}{N} \right) e \left( \frac{ajq_1^2}{q^2} \right) \\
\leq & \frac{\pi^2Q^3}{2N} + \frac{\pi^2 Q^2}{N} \sum_{0 < j < \frac{N}{4Q^2}} \left| \sum_{Q < q_1 \leq 2Q} e \left( \frac{ajq_1^2}{q^2} \right) \right|,
\end{eqnarray*}
where the first term above corresponds to the contribution of $j=0$.  Applying Cauchy-Schwartz inequality, we see that the square of the above expression is bounded by
\begin{equation*}
\ll \frac{Q^6}{N^2} + \frac{Q^2}{N} \sum_{0 < j < \frac{N}{4Q^2}} \left| \sum_{Q < q_1 \leq 2Q} e \left( \frac{ajq_1^2}{q^2} \right) \right|^2.
\end{equation*}
Applying Weyl Shift, Lemma~\ref{weylk} to the inner-most sum of the second term, we see that the double sum of the second term is
\begin{equation*}
\ll \sum_j Q + \sum_j \sum_{0 < l < Q} \min \left\{ Q, \left\| \frac{2ajl}{q^2} \right\|^{-1} \right\}
\ll \frac{N}{Q} + \sum_{0 < m < NQ^{-1}} \tau (m) \min \left\{ Q, \left\| \frac{am}{q^2} \right\|^{-1} \right\},
\end{equation*}
where $\tau(m)$ is the divisor function, is $O(m^{\epsilon})$ and estimates the multiplicity of representations of $m=2jl$.  The inequalities go in the correct direction by the virtue of positivity. \newline

What still remains is to estimate the sum over $m$.  We have, with $am \equiv d \pmod{q^2}$ and $|d| \leq \frac{1}{2} q^2$,
\begin{eqnarray*}
\sum_m & \leq & \sum_{|d| < q^2} \min \left\{ Q, \frac{q^2}{|d|} \right\} \sum_{0 < m < NQ^{-1}} \tau(m) \\
 & \ll & q^{-2} N^{1+\epsilon} + \sum_{0 < d < q^2} \frac{q^2}{d} \left( q^{-2} \frac{N}{Q} +1 \right) \left( \frac{N}{Q} \right)^{\epsilon} \\
 & \ll & q^{-2} N^{1+\epsilon} + \left( \frac{N}{Q} + q^2 \right) N^{\epsilon} \\
 & \ll & \left( \frac{N}{Q} + Q^2 \right) N^{\epsilon}.
\end{eqnarray*}
Recall that we are only considering the $q$'s in the dyadic interval $Q < q \leq 2Q$.  Combining everything and taking the square root, we infer that for every $x \in S_Q$,
\begin{equation*}
\# \left\{ x' \in S_Q : \| x-x' \| < \frac{1}{2N} \right\} \ll \frac{Q^3}{N} + \left( \sqrt{Q} + \frac{Q^2}{\sqrt{N}} \right) N^{\epsilon},
\end{equation*}
from which we infer the lemma.  \end{proof}

Now we are able to state and prove our main contention of the paper.  The beginning of the proof will go very much like that of the classical large sieve inequalities.  As far as that part is concerned, we are following the proof given in \cite{HI4}.

\begin{theorem} \label{larsiev}
With $\{ a_n \}$, $Q$, $M$, and $N$ defined as before, we have
\begin{equation} \label{larsieveq}
\sum_{q=1}^Q \sum_{\substack{a  \; \bmod{ \; q^2} \\ \gcd(a,q)=1}} \left| \sum_{n=M+1}^{M+N} a_n e \left( \frac{a}{q^2} n \right) \right|^2 \ll  \log 2Q \left[ Q^3+ (N\sqrt{Q} +\sqrt{N}Q^2)N^{\epsilon} \right] \sum_{n=M+1}^{M+N} |a_n|^2,
\end{equation}
where the implied constant depends on $\epsilon$ alone.
\end{theorem}

\begin{proof} It is easily observed that the theorem, after breaking the summation over $q$ into dyadic intervals with $Q < q \leq 2Q$ and the duality principle Lemma~\ref{dual}, and by assuming $M=0$ via the shift, $n \longrightarrow n-M$, it suffices to show that 
\begin{equation} \label{maincont}
\sum_{0 < n \leq N} \left| \sum_{x\in S_Q} b_x e(x n) \right|^2 \ll \left[ Q^3+ (N\sqrt{Q} +\sqrt{N}Q^2)N^{\epsilon} \right] \sum_{x\in S_Q} |b_x|^2,
\end{equation}
for any sequence of complex numbers $\{ b_x \}$. \newline

As before, we take $\phi(x) = \left( \frac{\sin \pi x}{2x} \right)^2$.  By positivity, the left-hand side of \eqref{maincont} is majorized by
\begin{equation*}
\sum_{n=-\infty}^{\infty} \phi \left( \frac{n}{2N} \right) \left| \sum_{x \in S_Q} b_x e (nx) \right|^2 = \sum_x \sum_{x'} b_x \overline{b}_{x'} V(x-x'),
\end{equation*}
where $V(y) = \sum_n \phi \left( \frac{n}{2N} \right) e(ny)$. \newline

Apply Poisson summations formula and a change of variables,
\begin{equation*}
V(y) = 2N \sum_m \hat{\phi} [2N (m+y)] = \frac{\pi^2N}{2} \sum_{|m+y|< (2N)^{-1} } ( 1-2N|m+y|) = \frac{\pi^2N}{2} (1-2N \|y\|),
\end{equation*}
if $\| y \| < (2N)^{-1}$ and $V(y)=0$ otherwise. \newline

Hence the left-hand side of \eqref{maincont} is majorized by
\begin{equation*}
\frac{\pi^2N}{2} \mathop{\sum_x \sum_{x'}}_{\|x-x'\|< (2N)^{-1}} |b_xb_{x'}| \leq \frac{\pi^2N}{2} \sum_x |b_x|^2 M(Q,N),
\end{equation*}
where $M(Q,N)$ is as defined in Lemma~\ref{lsscent}.  We insert the result of the afore-mentioned lemma, the theorem is proved.
\end{proof}

The greatest strength of our result lies in the range where $N \asymp Q^3$.  There, our result gives the majorant of $O(Q^{7/2+\epsilon})$ while both \eqref{Q4} and \eqref{Q3} give the majorant of $O(Q^4)$. \newline

In the same spirit that the classical large sieve inequality for additive characters implies that of the multiplicative characters, we have the following easy corollary.

\begin{corollary} \label{larsievcor}
For any sequence of complex numbers $\{ a_n \}$, we have
\begin{equation}
\sum_{q=1}^Q \frac{q}{\varphi (q)} \sideset{}{^{\star}}\sum_{\chi\; \bmod{ \, q^2}} \left| \sum_{n=M+1}^{M+N} a_n \chi (n) \right|^2 \ll \log 2Q \left[ Q^3 + (N\sqrt{Q}+\sqrt{N}Q^2)  N^{\epsilon} \right] \sum_{n=M+1}^{M+N} |a_n|^2,
\end{equation}
where the implied constant depends on $\epsilon$ alone.
\end{corollary}
\begin{proof} Using Gauss sums $G(\chi)$, we have
\begin{equation*}
\chi (n) = \frac{1}{G (\overline{\chi})} \sum_{a  \; \bmod{ \; q^2}} \overline{\chi} (a) e \left( \frac{an}{q^2} \right).
\end{equation*}
It is an elementary fact that the modulus of the Gauss sum $G (\chi)$ is the square root of its modulus, $q$ in our case.  Hence we have,
\begin{equation*}
\sideset{}{^{\star}}\sum_{\chi  \; \bmod{ \; q^2}} \left| \sum_{n=M+1}^{M+N} a_n \chi (n) \right|^2 \leq \frac{1}{q^2} \sum_{\chi \; \bmod{ \; q^2}} \left| \sum_{a \; \bmod{ \; q^2}} \overline{\chi} (a) \sum_n a_n e \left( \frac{an}{q^2} \right) \right|^2 = \frac{\varphi(q)}{q} \sum_{\substack{a \; \bmod{ \; q^2} \\
\gcd(a,q)=1}} \left| \sum_n a_n e \left( \frac{an}{q^2} \right) \right|^2.
\end{equation*}
The last equality is obtained by opening the modulus square and applying the
orthogonality of Dirichlet characters.  Note also that $\varphi(q^2) = q \varphi(q)$.  Our the corollary now follows from Theorem~\ref{larsiev}.
\end{proof}

We do not believe \eqref{larsieveq} is the best possible result.  In fact, we shall make conjectures based on computational evidences in Section 6, which yields a stronger result.

\section{Higher power moduli}

We mentioned that our method of the previous section may be generalized to investigate large sieve inequalities with higher power moduli.  Hence we dispose of that comment with the statement of our result here.  ``Since brevity is the soul of wit and tediousness the limbs and outward flourishes, I will be brief.''(\cite{WS}, {\it Hamlet, II, ii, 91-93}).  Hence we shall not provide all details of the proof as the proof goes {\it \`a la} previous section.  First we are writing down the lemma concerned the spacing of special fractions.

First set
\begin{equation*}
S_{Q,k} = \left\{ \frac{a}{q^k} \in \ratq : \gcd(a,q)=1, \; 1 \leq a < q^k, \; Q < q \leq 2Q \right\}, \; \mbox{and}
\end{equation*}
\begin{equation*}
M_k(Q,N) = \max_{x' \in S_{Q,k}} \# \left\{ x \in S_{Q,k} : \| x-x' \| < \frac{1}{2N} \right\}.
\end{equation*}
We have the following Lemma.

\begin{lemma} \label{lsscentk}
Given $\epsilon >0$, $N \in \natn$ and setting $\kappa = 2^{k-1}$, we have
\begin{equation} \label{lsscenteq2}
M_k(Q,N) \ll \frac{Q^{k+1}}{N} + \left( Q^{\frac{\kappa-1}{\kappa}}+\frac{Q^{\frac{\kappa+k}{\kappa}}}{N^{\frac{1}{\kappa}}} \right) N^{\epsilon},
\end{equation}
where the implied constant in \eqref{lsscenteq2} depends on $\epsilon$ and $k$.
\end{lemma}

From the above lemma, we have the following theorem.

\begin{theorem} \label{larsievk}
Let $\{ a_n \}$ be an arbitrary sequence of complex numbers, $Q$, $M$, $N$ positive integers, and $k \geq 2$.  Set $\kappa = 2^{k-1}$, we have
\begin{equation} \label{modk}
\sum_{q=1}^Q \sum_{\substack{a \; \bmod{ \; q^k} \\ \gcd(a,q)=1}} \left| \sum_{n=M+1}^{M+N} a_n e \left( \frac{a}{q^k} n \right) \right|^2 \ll \log 2Q \left[ Q^{k+1} + N^{\epsilon} \left( NQ^{\frac{\kappa-1}{\kappa}} + N^{1-\frac{1}{\kappa}} Q^{\frac{\kappa+k}{\kappa}} \right) \right] \sum_{n=M+1}^{M+N} |a_n|^2,
\end{equation}
where the implied constant will depend on $k$ and $\epsilon$.
\end{theorem}

Of course, in the same spirit as Corollary~\ref{larsievcor}, Theorem~\ref{larsievk} also admits the following corollary for multiplicative characters.

\begin{corollary}
With the same notation as before, we have
\begin{equation}
\sum_{q=1}^Q \frac{q}{\varphi (q)} \sideset{}{^{\star}}\sum_{\chi  \; \bmod{ \; q^k}} \left| \sum_{n=M+1}^{M+N} a_n \chi (n) \right|^2 \ll \log 2Q \left[ Q^{k+1} + N^{\epsilon} \left( NQ^{\frac{\kappa-1}{\kappa}} + N^{1-\frac{1}{\kappa}} Q^{\frac{\kappa+k}{\kappa}} \right) \right] \sum_{n=M+1}^{M+N} |a_n|^2,
\end{equation}
where the implied constant depends on $\epsilon$ and $k$.
\end{corollary}

As an application of the classical large sieve inequality, we can take
\begin{equation} \label{del2k}
\Delta = Q^{2k}+N,
\end{equation}
as all fractions with $k$-power denominator not greater than $Q^k$ will be $Q^{-2k}$-spaced modulo 1.  Alternatively, as before, since fractions with height exactly $q^k$ between 0 and 1, for $1\leq q \leq Q$ are $q^{-k}$-spaced, we may take
\begin{equation} \label{delk1}
\Delta = Q(Q^k+N).
\end{equation}
It is clear that the virtue of \eqref{del2k}, for large $k$'s, is limited when compared with \eqref{delk1}. \newline

As before, our result is most useful when $N \asymp Q^{k+1}$. In this case, \eqref{delk1} gives $\Delta=O(Q^{k+2})$, and \eqref{del2k} gives the extremely poor $\Delta=O(Q^{2k})$, while our result gives $\Delta=O(Q^{k+1+\frac{\kappa-1}{\kappa}})$, where $\kappa$ is defined as before.

\section{Conjectures based on heuristics and empirical evidences}

If we take $N=Q^3$ in Lemma~\ref{lsscent}, we have the following
\begin{equation} \label{lsscent2}
M(Q) \stackrel{\mathrm{def}}{=} \max_{x' \in S_Q} \# \left\{ x \in S_Q : 2 \| x-x' \| < Q^{-3} \right\} \ll Q^{\frac{1}{2}+\epsilon},
\end{equation}
and the implied constant in \eqref{lsscent2} depends on $\epsilon$ alone.  With this, we can arrive at something a slightly different.

\begin{prop} \label{larsiev2}
With $\{ a_n \}$, $Q$, $M$, and $N$ defined as before, we have
\begin{equation} \label{larsieveq2}
\sum_{q=1}^Q \sum_{\substack{a  \; \bmod{ \; q^2} \\ \gcd(a,q)=1}} \left| \sum_{n=M+1}^{M+N} a_n e \left( \frac{a}{q^2} n \right) \right|^2 \ll  Q^{\frac{1}{2} + \epsilon} (Q^3+ N) \sum_{n=M+1}^{M+N} |a_n|^2,
\end{equation}
where the implied constant depends on $\epsilon$ alone.
\end{prop}

\begin{proof} The proposition follows easily from Theorem~\ref{larsiev} with the observation that $Q^2 \sqrt{N} \leq \max \{Q^{\frac{7}{2}}, \sqrt{Q} N \}$.  However, if one follows a different line of proof and partition the interval $(0,1]$ into $2Q^3$ subintervals of equal size, say $K_i, \; i = 1, \cdots , 2Q^3$, then $\# \{K_i \cap S_Q\} = O(Q^{\frac{1}{2}+\epsilon})$, for all $i$, by the virtue of \eqref{lsscent2}.  Picking one element from $K_i \cap S_Q$ for each even and then odd $i$ would ensure the well-spacedness of the resulting set and such a ``picking'' process terminates after $O(Q^{\frac{1}{2}+\epsilon})$ times.  Thus, one can also infer the result of the proposition and the factor of $Q^{\frac{1}{2}+\epsilon}$ \eqref{larsieveq2} comes from that in \eqref{lsscent2}. \end{proof}

The above gives the significance of the majorant \eqref{lsscent2} which will facilitate our conjecture making.  Table 1 below lists some values of $M(Q)$ for some small values of $Q$, where $M(Q)$ is as in \eqref{lsscent2}.  We see that $M(Q)$ increases rather slowly and it is our belief that $M(Q)=O(Q^{\epsilon})$ with the implied constant depending on $\epsilon$ alone.  Indeed, there are about $Q^3$ elements in $S_Q$, and ``on average,'' they are $Q^{-3}$-spaced.  Hence, there should not be many pairs of elements in $S_Q$ that are spaced too closely.  The readers should allow my making of the following analogue.  If we are using a sieve with $Q^3$ holes of equal size in it and sift out the elements of $S_Q$ with it, then ``not many'' rational numbers should fall through the same hole.  This is precisely the meaning of the conjectures in this section. \newline

\begin{table} \label{tab1} \caption{$M(Q)$ for some small $Q's$.}
\begin{center}
\begin{tabular}{|c|c|c|c|c|c|c|c|c|c|}\hline 
$Q$ & $M(Q)$ & $Q$ & $M(Q)$ & $Q$ & $M(Q)$ & $Q$ & $M(Q)$ & $Q$ & $M(Q)$ \\ \hline
1 & 0 & 2 & 0 & 3 & 1 & 4 & 1 & 5 & 2 \\ \hline
6 & 1 & 7 & 1 & 8 & 2 & 9 & 2 & 10 & 2 \\ \hline
11 & 2 & 12 & 2 & 13 & 2 & 14 & 2 & 15 & 2 \\ \hline
16 & 2 & 17 & 2 & 18 & 2 & 19 & 2 & 20 & 2 \\ \hline
21 & 2 & 22 & 2 & 23 & 2 & 24 & 3 & 25 & 2 \\ \hline
26 & 2 & 27 & 3 & 28 & 2 & 29 & 2 & 30 & 2 \\ \hline
31 & 2 & 32 & 2 & 33 & 2 & 34 & 2 & 35 & 2 \\ \hline
36 & 2 & 37 & 3 & 38 & 3 & 39 & 3 & 40 & 3 \\ \hline
41 & 3 & 42 & 3 & 43 & 3 & 44 & 3 & 45 & 2 \\ \hline
46 & 3 & 47 & 3 & 48 & 3 & 49 & 3 & 50 & 3 \\ \hline
51 & 3 & 52 & 3 & 53 & 4 & 54 & 3 & 55 & 3 \\ \hline
56 & 3 & 57 & 3 & 58 & 3 & 59 & 3 & 60 & 3 \\ \hline
61 & 3 & 62 & 3 & 63 & 3 & 64 & 3 & 65 & 3 \\ \hline
66 & 3 & 67 & 3 & 68 & 3 & 69 & 3 & 70 & 3 \\ \hline
71 & 3 & 72 & 3 & 73 & 3 & 74 & 3 & 75 & 3 \\ \hline
76 & 3 & 77 & 3 & 78 & 3 & 79 & 3 & 80 & 3 \\ \hline
81 & 3 & 82 & 3 & 83 & 4 & 84 & 4 & 85 & 3 \\ \hline
86 & 3 & 87 & 3 & 88 & 3 & 89 & 4 & 90 & 4 \\ \hline
91 & 4 & 92 & 4 & 93 & 4 & 94 & 3 & 95 & 3 \\ \hline
96 & 3 & 97 & 4 & 98 & 4 & 99 & 4 & 100 & 4 \\ \hline
\end{tabular}
\end{center}
\end{table}

We believe the growth of $M(Q)$, as defined in \eqref{lsscent2}, is of independent interest.  Therefore, we record the following conjectures.

\begin{conjecture} \label{conj1}
Let $S_Q$ be defined as before.  Then we have
\begin{equation}
\max_{x\in S_Q} \# \left\{ x' \in S_Q : \|x-x'\| < Q^{-3} \right\} \ll Q^{\epsilon},
\end{equation}
where the implied constant depends on $\epsilon$ alone.
\end{conjecture}

In the same spirit, we also express, but with less confidence, the following conjecture for higher power moduli.

\begin{conjecture} \label{conj2}
Let 
\[ S_{k,Q}= \left\{ \frac{a}{q^k} \in \ratq : \gcd(a,q)=1, \; 1 \leq a < q^k, \; Q < q \leq 2Q \right\}. \]
Then we have
\begin{equation}
\max_{x\in S_{k,Q}} \# \left\{ x' \in S_{k,Q} : \|x-x'\| < Q^{-k-1} \right\} \ll Q^{\epsilon},
\end{equation}
where the implied constant depends only on $k$ and $\epsilon$.
\end{conjecture}

If we assume the truth of Conjecture~\ref{conj2}, then we would have the following.
\begin{conjecture} \label{conj3}
Let $\{ a_n \}$ be an arbitrary sequence of complex numbers, $Q$, $N \in \natn$ and $M \in \intz$.  We have
\begin{equation} \label{conj3eq}
\sum_{q=1}^Q \sum_{\substack{a=1 \\ \gcd(a,q)=1}}^{q^k} \left| \sum_{n=M+1}^{M+N} a_n e \left( \frac{a}{q^k} n \right) \right|^2 \ll Q^{\epsilon} (Q^{k+1}+ N) \sum_{n=M+1}^{M+N} |a_n|^2,
\end{equation}
where the implied constant depends only on $\epsilon$ and $k$.
\end{conjecture}

\bibliography{biblio}
\bibliographystyle{amsxport}

\vspace*{.1in}
\hspace*{.5in}{\sc\small Dept. Math., Rutgers Univ., 100 Frelinghuysen Rd., Piscataway, NJ 08854 USA \newline
\hspace*{.5in}\indent Email Address: {\tt lzhao@math.rutgers.edu}

\end{document}